\documentclass[10pt]{amsart}
\usepackage{graphics}
\usepackage{amsfonts}
\usepackage{amsmath}

\title{El polinomio de Jones y la mec\'{a}nica cu\'{a}ntica}
\author{R{\u{a}}zvan Gelca}
\address{Department of Mathematics and Statistics, 
Texas Tech University, Lubbock, TX 79409 and Institute of Mathematics
of the Romanian Academy, Bucharest, Romania}
\email{rgelca@gmail.com}

\subjclass{81S10, 81R50, 57R56, 81T45, 57M25}
\keywords{Jones polynomial, topological quantum
field theory, 
moduli spaces of flat connections, quantization}

\begin{document}

\begin{abstract}
From the moment of its discovery, the Jones polynomial of a knot has been linked
to quantum physics. The main discovery, made  by E. Witten,
was that it is related to quantum field theory, which unfortunately lacks
a mathematical foundation. But already in Witten's work it was noted that the
Jones polynomial is related to quantum mechanics. 

In this paper we discuss progress made in the study of the Jones polynomial
from the point of view of quantum mechanics. This study
 reduces to the understanding
of the quantization of the moduli space of flat $SU(2)$-connections on a 
surface with  the Chern-Simons lagrangian. We outline some background material, then
present the particular example of the torus, in which case it is known that
the quantization in question is the Weyl quantization. The paper concludes with
a possible application of this theory to the study of the
fractional  quantum Hall effect,  an idea originating in the works of
Moore and Read. 
\end{abstract}

\maketitle

\section{Resumen}

En este  art\'{\i}culo de exposici\'{o}n  se presenta una introducci\'{o}n al 
 estudio del polinomio de Jones 
de un nudo desde un punto de
vista f\'{\i}sico basado en la 
{\em mec\'{a}nica cu\'{a}ntica}.\footnote{\'{E}ste art\'{\i}culo ha sido
publicado en Aportaciones Mat. Comun., 36, Soc. Mat. Mexicana, 2006, 85--99.} 

El polinomio de Jones de un nudo fue descubierto por V.F.R. Jones
\cite{jones} como una consecuencia de sus investigaciones en la
teor\'{\i}a de \'{a}lgebras de operadores. La  definici\'{o}n
del polinomio es muy sencilla y nos permite calcularlo facilmente, pero
para entender y estudiar sus propiedades se necesita un punto 
de vista mas profundo y mas geom\'{e}trico.
 Witten explic\'{o}  \cite{witten} que 
el polinomio de Jones esta relacionado a la teor\'{\i}a de campos
cu\'{a}nticos.
Pero la teor\'{\i}a de campos cu\'{a}nticos requerida no tiene un 
 fundamento matem\'{a}tico riguroso al momento de escribir esto.   

En su art\'{\i}culo, Witten explic\'{o} que el estudio del polinomio
de Jones se puede reducir a un problema de mec\'{a}nica cu\'{a}ntica.
El problema es de cuantizar el espacio de 
campos con simetr\'{\i}a interior $SU(2)$ en una superficie modulo
cambios locales de coordenadas, que es el mismo que el espacio de moduli
de conexiones planas de $su(2)$ en la superficie.
 El espacio de Hilbert
de esta cuantizaci\'{o}n se entiende bien, y consiste de las secciones
holom\'{o}rficas del haz en l\'{\i}neas de Chern-Simons (que se llaman 
funciones theta generalizadas) \cite{freed}. 
Por otro lado, aunque hay varios trabajos en esta direcci\'{o}n
\cite{alexeevschomerus} \cite{jea}, \cite{GU}, 
los operadores de la cuantizaci\'{o}n (los observables
cu\'{a}nticos) son un misterio.
Todav\'{\i}a no hay una respuesta
a la pregunta principal: ?`Que procedimiento de 
cuantizaci\'{o}n corresponde al polinomio de Jones? 

En lo que sigue presentaremos la definici\'{o}n del polinomio
de Jones y su relaci\'{o}n con la teor\'{\i}a de campos cu\'{a}nticos.
Luego explicaremos la reducci\'{o}n de esta teor\'{\i}a 
de campos cu\'{a}nticos al problema de cuantizaci\'{o}n del espacio
de moduli de conexiones planas de $su(2)$ en una supeficie 
y revisamos las cuestiones b\'{a}sicas de  mec\'{a}nica cu\'{a}ntica. 
El art\'{\i}culo sigue con una secci\'{o}n que 
discute la estructura del espacio de moduli de conexiones planas de $su(2)$  
 en una superficie 
y presenta el progreso
realizado para entender su cuantizaci\'{o}n. 
Como caso particular  discutiremos el ejemplo del toro que es de los m\'{a}s estudiados,
para el cual presentamos las relaciones con el toro no conmutativo
y con la cuantizaci\'{o}n de H. Weyl.
El art\'{\i}culo concluye con
 una posible aplicaci\'{o}n de la teor\'{\i}a f\'{\i}sica del
polinomio de Jones 
a la teor\'{\i}a del efecto  Hall cu\'{a}ntico fraccionario.

El autor quiere agradecer a Ernesto Lupercioy a Lourdes Juan por ayudalo mejorar el estilo de esta exposici\'{o}n.

\section{El polinomio de Jones y la teor\'{\i}a de los campos 
cu\'{a}nticos}

El polinomio de Jones es un  
invariante de nudos descubierto por
   V.F.R. Jones \cite{jones}. Este invariante esta definido por la relaci\'{o}n
recurrente 
\begin{eqnarray*}
{t^{-1}V_{K_+}(t)-tV_{K_-}(t)=(t^{1/2}-t^{-1/2})V_{K_0}(t)},
\end{eqnarray*}
en la cual los diagramas de los 
 nudos $K_+$, $K_-$ y $K_0$ son los mismos salvo por un cruzamiento
como es represantado en la Figura 1. Esta relaci\'{o}n  nos permite
calcular  el polinomio de Jones $V_K(t)$ de cualquier  nudo $K$ sabiendo 
el  polinomio  de nudos mas sencillos. Cambiando cruzamientos como
se muestra cada nudo se puede transformar en el nudo trivial, cuyo polinomio
es $V_0(t)=1$.
La Figura 2 ejemplifica el c\'{a}lculo del polynomio de Jones para
el nudo tr\'{e}bol. Al final de este c\'{a}lculo se obtiene
$V_{\mbox{tr\'{e}bol}}(t)=t+t^3-t^4$.

\begin{figure}[htbp]
\centering
\scalebox{.30}{\includegraphics{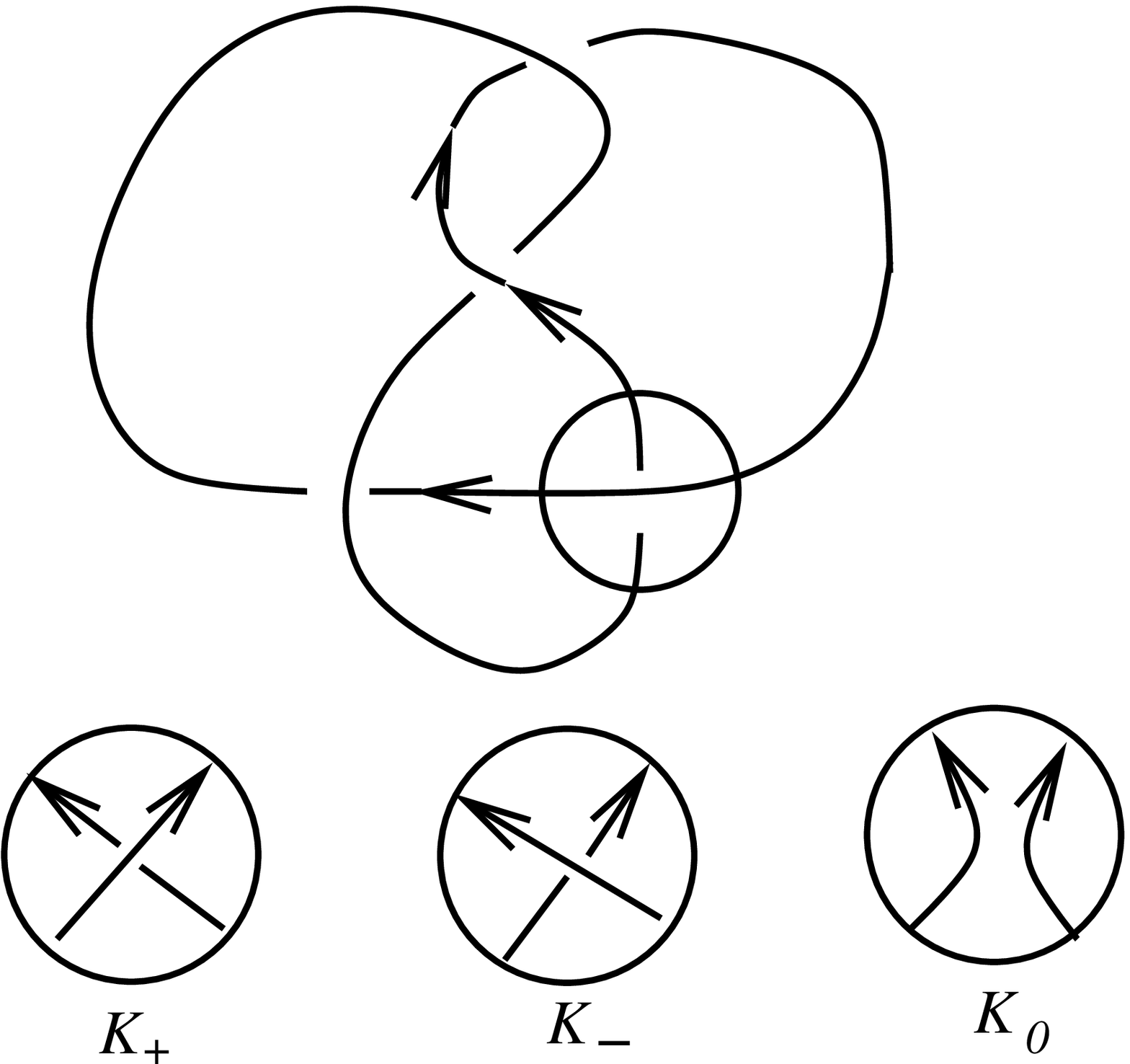}}

Figura 1.
\end{figure}

\begin{figure}[h]
\centering
\scalebox{.30}{\includegraphics{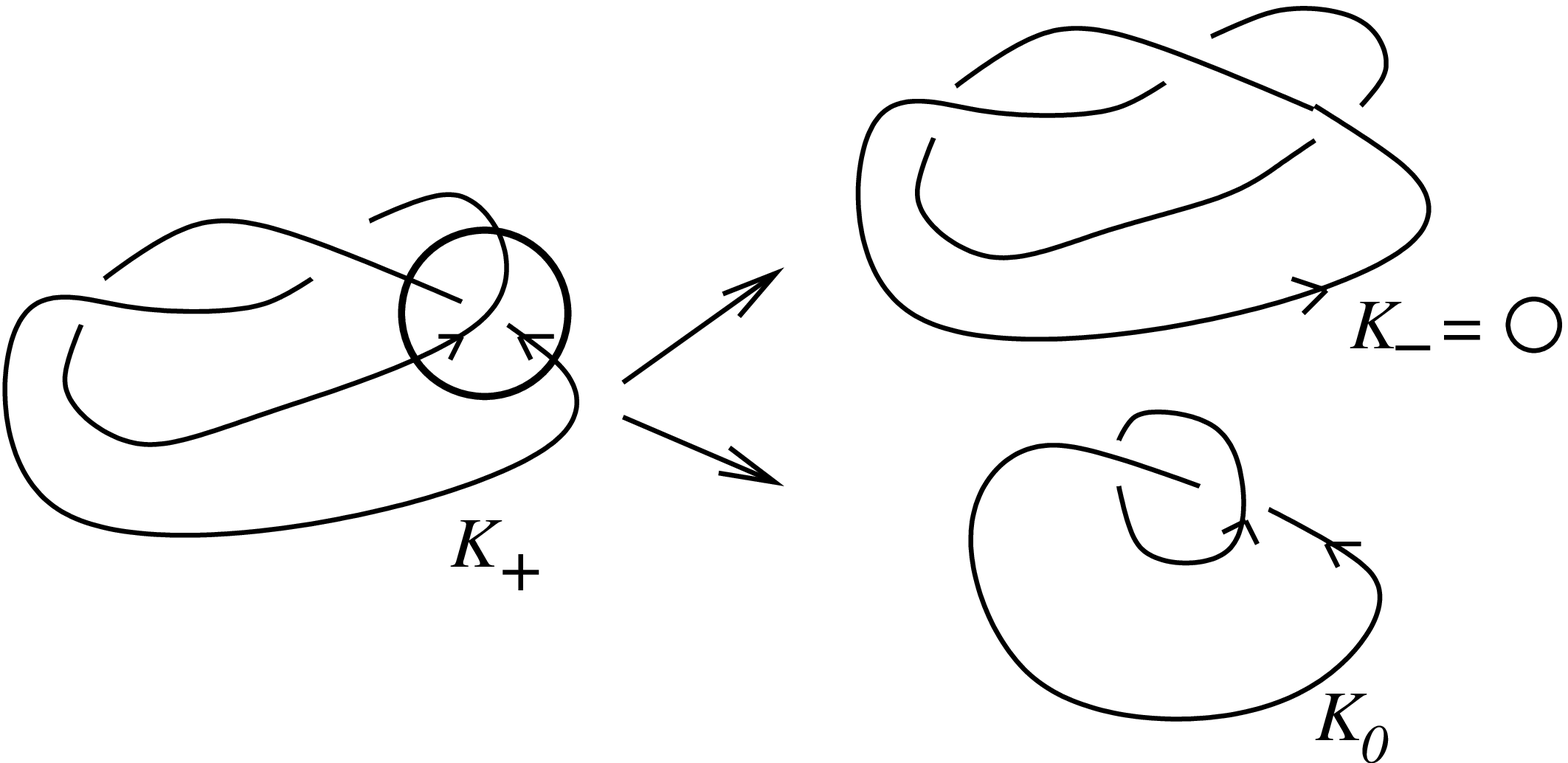}}\\
$ V_{K_+}(t)=t^2+(t^{3/2}-t^{1/2})V_{K_0}(t)$\\
\scalebox{.30}{\includegraphics{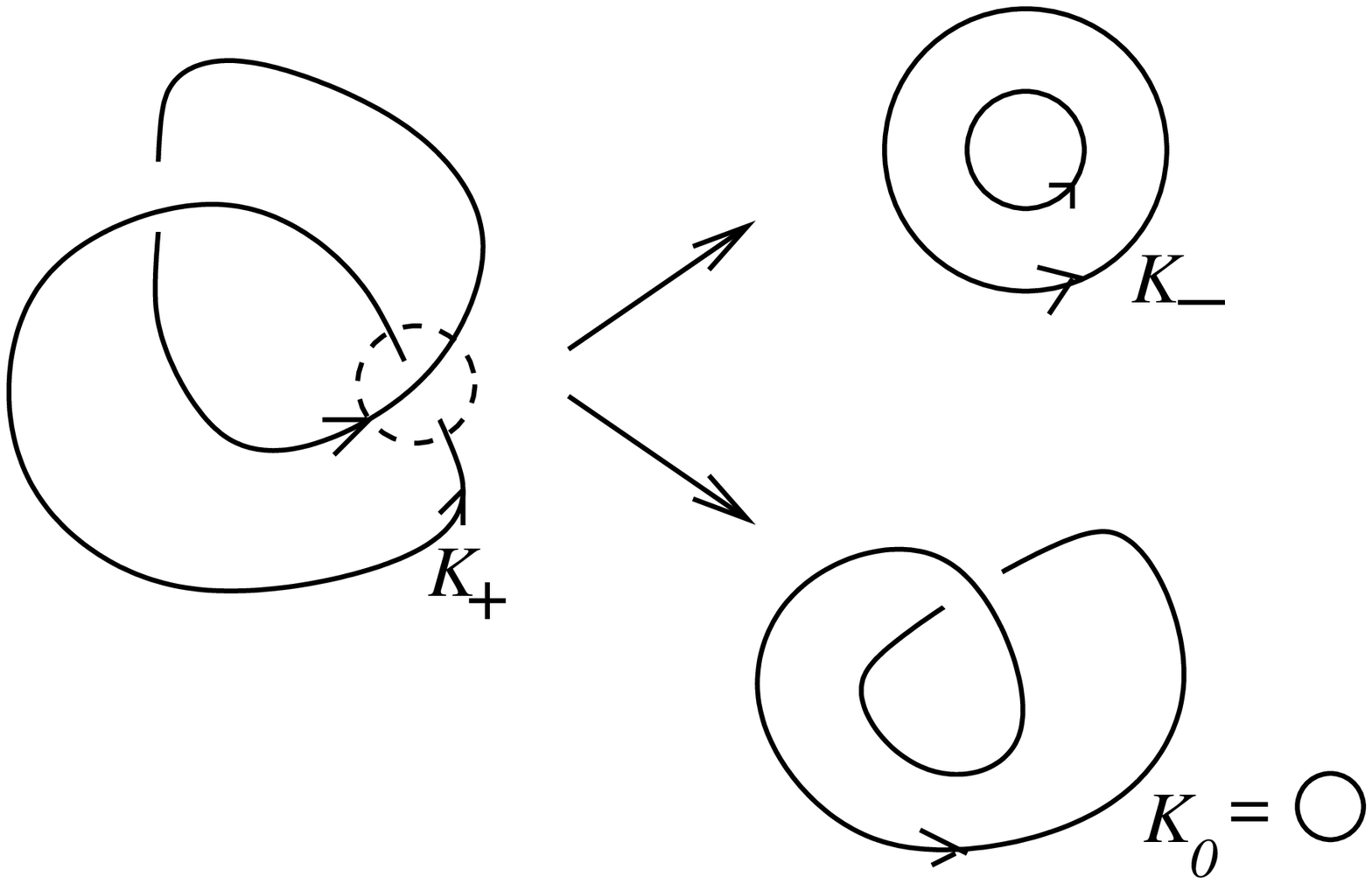}}
\scalebox{.30}{\includegraphics{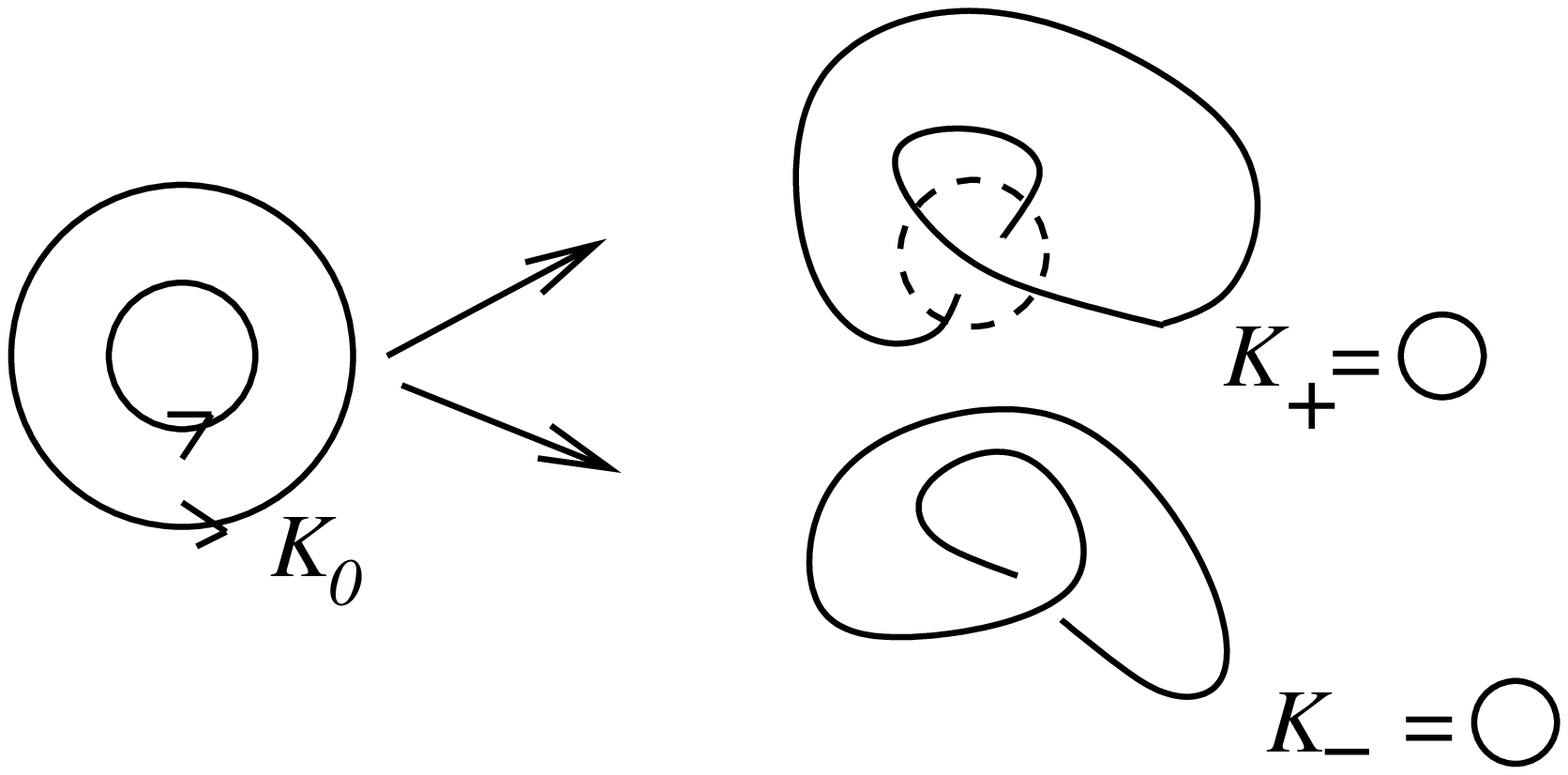}}\\
$V_{K_+}(t)=t^2V_{K_-}(t)+(t^{3/2}-t^{1/2}),\quad V_{K_0}(t)=t^{1/2}+t^{-1/2}.$\\ 

\medskip

Figura 2.
\end{figure}

Esta definici\'{o}n es muy \'{u}til para calcular el  polinomio, pero no
permite el estudio de sus propiedades. Los top\'{o}logos
desean una definici\'{o}n mas geom\'{e}trica, de cual se pueden deducir
algunas de las propriedades geom\'{e}tri\-cas y topol\'{o}gicas del nudo. 

Cinco a\~{n}os despues de Jones, 
E. Witten \cite{witten} descubri\'{o}
que el polinomio de Jones  aparece en la
teor\'{\i}a de campos cu\'{a}nticos. Witten tuvo la idea siguente.

Se consideran los campos con grupo de simetr\'{\i}a interna
\begin{eqnarray*}
{SU(2)}= \left\{
\left[
\begin{array}{rr}
a&b\\
-\bar{b} &\bar{a}
\end{array}
\right],\, a,b\in {\mathbb C}, 
a\bar{a}+b\bar{b}=1\right\}
\end{eqnarray*}
en una variedad $M$ de dimension tres.
Cada campo tiene un potencial,
que desde el punto de vista matem\'{a}tico es una  
 conexi\'{o}n $A$ de $su(2)$ en $M$.
La presencia del campo se puede determinar por su acci\'{o}n en una
 part\'{\i}cula que se mueve en el espacio.
La fase de la part\'{\i}cula, que es un vector en ${\mathbb R}^2$
(ver la Figura 3), que cambia de acuerdo con
\begin{eqnarray*}
\left[
\begin{array}{r}
x\\
y
\end{array}
\right]
\rightarrow 
{\left[
\begin{array}{rr}
a&b\\
-\bar{b} &\bar{a}
\end{array}
\right]}
\left[
\begin{array}{r}
x\\
y
\end{array}
\right]. 
\end{eqnarray*}
La matriz de rotaci\'{o}n se llama   {matriz de holonom\'{\i}a}.
\begin{figure}[htbp]
\centering
\scalebox{.30}{\includegraphics{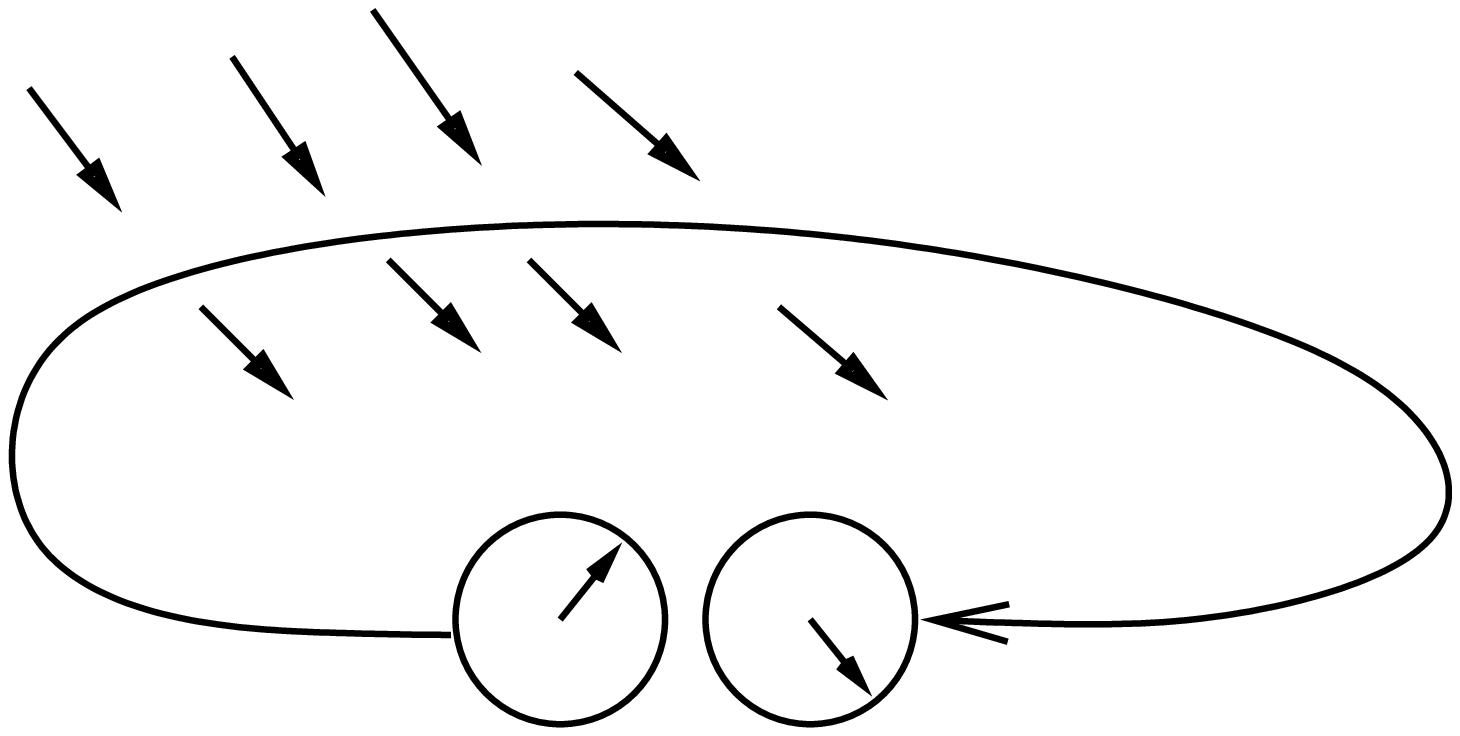}}

Figura 3.
\end{figure}

El lagrangiano de la teor\'{\i}a de campo es el lagrangiano de
Chern-Simons
 \begin{eqnarray*}
{L(A)}=
\frac{1}{4\pi}\int_{M}tr(A\wedge dA+\frac{2}{3}A\wedge A\wedge
A).
\end{eqnarray*}
 
De los observables f\'{\i}sicos, nos interesa la
traza de la holonom\'{\i}a, que se denota ${W_K(A)}$.
Como observable cl\'{a}sico este es una funci\'{o}n de
la conexi\'{o}n, como observable cu\'{a}ntico este
es un operador lineal. Para cuantizar esta teor\'{\i}a de campos introducimos
la constante  de Planck $\hbar$ que satisface la condici\'{o}n
de cuantizaci\'{o}n geometrica $\hbar=1/N$ con $N$ un entero positivo,
y adem\'{a}s pedimos que $N$ sea par: $N=2r$. 

\begin{figure}[htbp]
\centering
\scalebox{.40}{\includegraphics{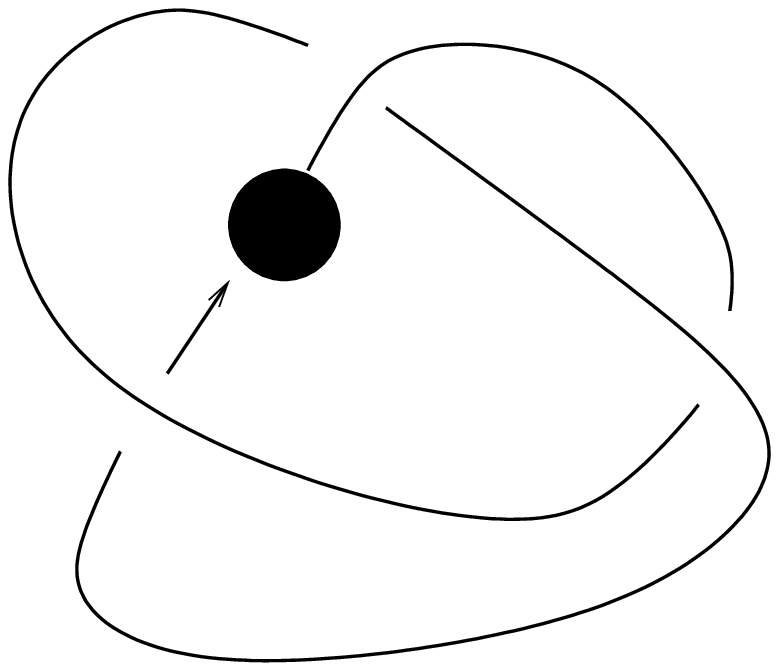}}

Figura 4.
\end{figure}

Para una part\'{\i}cula que se mueve por  un 
nudo como el mostrado en la Figura 4, se toma la
media de la traza de  la holonom\'{\i}a 
sobre todos los campos. Si la variedad $M$ es la esfera de dimension
tres, el resultado es  el  polinomio de Jones del  nudo evaluado
en una ra\'{\i}z $r$-\'{e}sima de la unidad.
Formalmente, esta media es la integral de Feynman 
\begin{eqnarray*}
\int e^{2i r{L(A)}}{W_K(A)}{\mathcal D}A.
\end{eqnarray*}
Porque ambas cantidades  $e^{2i r{L(A)}}$ y $W_K(A)$ se mantienen invariantes
con el cambio de sistema de coordenadas, la media se toma sobre
clases de campos equivalentes por cambios de coordenadas.
En conclusi\'{o}n, el polinomio de Jones es un objecto de la
teor\'{\i}a de campos cu\'{a}nticos.

Desafortunadamente, la teor\'{\i}a general de campos cu\'{a}nticos no tiene
un fundamento matem\'{a}tico, pues la integral de Feynman no se puede
definir desde un punto de vista matem\'{a}tico. 
El problema es que un campo tiene una infinitud de 
grados de libertad,
y  en este momento comprendemos solamente quantizaciones de sistemas
con un numero finito de grados de libertad.   

Para el polinomio de Jones hay una teor\'{\i}a desarrollada por Reshetikhin y 
Turaev \cite{RT} que imita la teoria de Witten,
pero que pierde el sabor f\'{\i}sico y geom\'{e}trico. 
Su teor\'{\i}a usa los grupos cu\'{a}nticos de Drinfeld \cite{drinfeld},
con c\'{a}lculos combinatorios muy complicados.   

Por otro lado, en nuestro caso los campos se identifican por 
cambios  locales de coordenadas (un procedimiento llamado \emph{gauging} o \emph{normalizaci\'{o}n}), 
obteniendo los {\em campos f\'{\i}sicos} en la variedad. 
Ademas, la teor\'{\i}a en 3 dimensiones se puede reducir
a una teor\'{\i}a en 2 dimensiones usando la descomposici\'{o}n
Heegaard de variedades de dimension tres.
Afortunadamente los campos f\'{\i}sicos en  una
 superficie tienen un numero {\em finito} de grados de libertad,
pues se pueden cuantizar con t\'{e}cnicas bien entendidas de 
mec\'{a}nica cu\'{a}ntica.
Esto nos entrega un problema de {mec\'{a}nica cu\'{a}ntica}:
hay que cu\'{a}ntizar el espacio de los campos con simetria 
$SU(2)$ en una superficie modulo cambios locales de coordenadas. 
La pregunta principal es,
de los varios modelos de 
cu\'{a}ntizacion cual corresponde al polinomio de Jones.

\section{Que es cu{a}ntizaci\'{o}n}

En esta secci\'{o}n explicamos algunas cuestiones b\'{a}sicas de 
 mec\'{a}nica cl\'{a}sica y cu\'{a}ntica. 
En la formulaci\'{o}n Hamiltoniana de la mec\'{a}nica cl\'{a}sica,
cada sistema tiene un
espacio  fase, que es el haz cotangente $T^*N$ de la variedad 
de configuraciones $N$,
o mas general, una variedad simpl\'{e}ctica $(M,\omega)$. 
Los observables f\'{\i}sicos son funciones $f\in C^{\infty}(M)$.
Por ejemplo si $N={\mathbb R}$, entonces $M=T^*N={\mathbb R}^2$ 
con coordenadas $q$ (posici\'{o}n) y $p$ (momento) y forma simpl\'{e}ctica
$dp\wedge dq$. Los observables son funciones $f(p,q)$, como por ejemplo la
energ\'{\i}a total $H=\frac{1}{2}p^2+\frac{1}{2}q^2$.  

La forma simpl\'{e}ctica asocia a cada observable $f$ un campo
vectorial $X_f$ definido por la ecuaci\'{o}n $df =\omega(X_f,\cdot)$.
De esta manera la forma simpl\'{e}ctica define un par\'{e}ntesis
de Poisson para funciones en $C^\infty(M)$  por la f\'{o}rmula
$\{f,g\}=df(X_g)$. Cada sistema f\'{\i}sico tiene
una funci\'{o}n de energ\'{\i}a total $H$, que junto con el
par\'{e}ntesis de Poisson, nos da la 
evoluci\'{o}n de un observable por medio de la ecuaci\'{o}n de Hamilton 
\begin{eqnarray*}
{\frac{df}{dt}=\{f,H\}}.
\end{eqnarray*}
El par\'{e}ntesis de Poisson para nuestro ejemplo particular 
es $\{f,g\}=\frac{\partial f}{\partial p}
\frac{\partial g}{\partial q}-\frac{\partial f}{\partial
  q}\frac{\partial g}{\partial p}$. Escribiendo la evoluci\'{o}n de
la posici\'{o}n y del momento se obtienen las ecuaciones cl\'{a}sicas
de Hamilton.

En la mec\'{a}nica cu\'{a}ntica de Heisenberg, los estados de un sistema
son los vectores  unitarios  en un 
espacio de Hilbert (funci\'{o}nes de onda). Los 
observables son operadores lineales 
Herm\'{\i}ticos actuando sobre el espacio de  Hilbert.
La observacion de un operador  ${T}$ en el estado $\psi$ 
no es determinista, su valor esperado  es
$|\left<{T}{\psi}, {\psi}\right>|$.

La evoluci\'{o}n de un observable $A$ es definida por medio de la
ecuaci\'{o}n de Schr\"{o}dinger 
\begin{eqnarray*}
{\frac{dA}{dt}=\frac{1}{i\hbar}[A,\widehat{H}]}=\frac{1}{i\hbar}(
A\widehat{H}-\widehat{H}A),
\end{eqnarray*}
donde $\widehat{H}$ es el operador de energ\'{\i}a total.

Cuantizaci\'{o}n es un procedimiento para pasar de la mec\'{a}nica 
cl\'{a}sica a la mec\'{a}nica cu\'{a}ntica.  
Este procedimiento substituye el espacio de fase con un espacio
de Hilbert (de funciones de onda) y las funciones sobre el espacio de
fase con  operadores actuando en el espacio de Hilbert.
La ecuaci\'{o}n de Hamilton se transforma en la ecuaci\'{o}n de 
Schr\"{o}dinger, con el par\'{e}ntesis de Poisson convertiendose en
el conmutador de operadores. 

En el caso m\'{a}s sencillo  del espacio de fase 
$T^*{\mathbb R}$, el espacio de Hilbert asociado por la cuantizaci\'{o}n
es  ${{\mathcal H}=L^2({\mathbb R})}$.
 En el caso general de la variedad cotangente $T^*(N)$, o el m\'{a}s general 
de la variedad simpl\'{e}ctica $(M, \omega)$ con $\omega$ una clase
de $H^2(M,{\mathbb Z})$, 
la cu{a}ntizaci\'{o}n geom\'{e}trica (ver \cite{woodhouse}) produce el espacio
de Hilbert ${\mathcal H}$ como el espacio $\Gamma (M, {\mathcal L})$ de las
secciones $L^2$-integrables de un haz complejo de l\'{\i}neas ${\mathcal L}$
sobre $M$.  
La forma simpl\'{e}ctica que describe la evoluci\'{o}n en el sistema
original y la constante de Planck $\hbar=\frac{1}{N}$
se incluyen en la estructura geom\'{e}trica del haz, pidiendo que
su curvatura sea $N\omega$. Varios casos particulares,
bien entendidos basado en  experimentos, indican que el
espacio de Hilbert definido as\'{\i} es
demasiado grande, por lo que lo reducimos   por
medio de una polarizaci\'{o}n. 

Por ejemplo, si identificamos $T^*{\mathbb R}$ con ${\mathbb C}$
poniendo  $z=q+ip$, el espacio de Hilbert
en la polarizaci\'{o}n real es el espacio de funciones $L^2$-integrables
en la variable real $q$. Si la variedad simpl\'{e}ctica
$M$ es una variedad compleja, hay lo que se llama la polarizaci\'{o}n 
holom\'{o}rfa, para cual el espacio de Hilbert
consiste de  las secciones holom\'{o}rfas $L^2$-integrables
del haz ${\mathcal L}$.

El producto interior de dos  secciones del espacio de Hilbert
 se puede definir como ${
\left<\sigma,\sigma'\right>}=\int_M \left<\sigma,
\sigma'\right>_x Vol$ donde $Vol=\omega \wedge\omega \wedge \cdots 
\wedge \omega$, pero es posible que las secciones producidas
por la polarizaci\'{o}n no sean integrables. El procedimiento est\'{a}ndar
consiste en introducir una semi-densidad relacionada a la polarizaci\'{o}n y
 multiplicar cada secci\'{o}n por esta semi-densidad. Desde el punto
de vista geom\'{e}trico, se multiplica el haz ${\mathcal L}$ por 
otro haz que contiene la informaci\'{o}n sobre la medida de integraci\'{o}n
(ver \cite{snyaticki}).

Los operadores lineales $op(f)$ asociados a funciones $f\in C^{\infty}(M)$
deben  respetar las  condiciones de Dirac:
\begin{itemize}
\item[1.] $op(1)=I$,
\item[2.] si $f$ es de valores reales, el operador $op(f)$ debe
ser autoadjunto.
\item[3.] $op(\{f,g\}) =\frac{1}{i\hbar}[op (f), op(g)]$.
\end{itemize}
En realidad la \'{u}ltima condici\'{o}n es demasiado optimista
y es imposible de satisfacer (un resultado de Groenewold y Van Hove), 
y por eso hay que introducir un  error $O(\hbar)$. 

En el caso m\'{a}s sencillo $T^*{\mathbb R}$ se requiere adem\'{a}s que
el operador $op(q)$ sea la multiplicaci\'{o}n por la coordenada $q$ y
el operador $op(p)$ sea $-i\hbar \frac{\partial}{\partial q}$, y que
el espacio de Hilbert sea irreducible bajo la acci\'{o}n de estos dos
operadores.
En este caso para definir el operador asociado a una funci\'{o}n general 
$f(p,q)$ se puede usar la ``receta'' de H. Weyl:
\begin{eqnarray*}
op(f)=\frac{1}{2\pi}\int_{{\mathbb R}^2}
\hat{f}(x,y)\exp(i(x\, op(p)+ y\, op
    (q)))dxdy,
\end{eqnarray*} 
donde $\hat{f}$ es la transformada de Fourier de $f$.
 Pero hay una infinitud de ``recetas'' de cuantizaci\'{o}n, como
por ejemplo la de Toeplitz:
\begin{eqnarray*}
op(f):\quad g\in L^2({\mathbb R})\longrightarrow fg\in L^2({\mathbb R}^2)
\longrightarrow \Pi fg\in  L^2({\mathbb R}),
\end{eqnarray*}
donde $\Pi$ es la proyecci\'{o}n de $L^2({\mathbb R}^2)$ en
el espacio de Hilbert $L^2({\mathbb R})$ de las funciones en
la primera variable.

El procedimiento de cuantizaci\'{o}n de Toeplitz se puede
generalizar  al caso general de una variedad simpl\'{e}ctica $M$: 
\begin{eqnarray*}
op(f): \quad \sigma\in {\mathcal H}\longrightarrow f\sigma \in L^2(M, L)
\longrightarrow \Pi f\sigma \in {\mathcal H}.
\end{eqnarray*}
Existe tambi\'{e}n el procedimiento de cuantizaci\'{o}n geometrica estandar
\cite{woodhouse}.
Otro metodo, desarollado en los ultimos 30 a\~{n}os \cite{bayenetall},
 es deformar el algebra de funciones $C^\infty(M)$
a un algebra no conmutativa de operadores:
\begin{eqnarray*}
f*g=fg +\hbar \{f,g\}+\hbar^2B_2(f,g)+\hbar^3B_3(f,g)+\cdots 
\end{eqnarray*}

\section{El espacio de moduli de conexiones planas y su cuantizaci\'{o}n}

Nuestra intenci\'{o}n es  cuantizar  los campos con simetr\'{\i}a interior $SU(2)$
en una superficie $S$ modulo cambios locales de coordenadas.
Desde el punto de vista matem\'{a}tico, el espacio de  fase es
\begin{eqnarray*}
\{A\, |\, A :\, su(2)-\mbox{conecci\'{o}n}\}/{\mathcal G}
\end{eqnarray*}
donde ${\mathcal G}$ es  el grupo de transformaciones gauge
\begin{eqnarray*}
A\rightarrow g^{-1}Ag+g^{-1}dg \quad \mbox{donde } g:S\rightarrow SU(2).
\end{eqnarray*}
Por  reducci\'{o}n simpl\'{e}ctica (ver \cite{atiyah})
 cuantizar este espacio es lo mismo que
cuantizar
\begin{eqnarray*}
{{\mathcal M}}{=\{A\, |\, A :\,
  su(2)-\mbox{conexi\'{o}n plana}\}/{\mathcal G}}
{=\{\rho: \pi_1(S)\longrightarrow SU(2)\}
/\mbox{conjugaci\'{o}n}}
\end{eqnarray*}
que es por un lado el espacio de moduli de conexiones planas de $su(2)$ en la
  superficie
y por otro lado la variedad de caracteres de $SU(2)$-representaciones del
grupo fundamental de la superficie.
${\mathcal M}$ es una variedad algebraica compleja muy complicada,
estudiada por primera vez por  
Narasimhan y Seshadri \cite{NS}.

Atiyah y Bott han mostrado \cite{AB} que
la parte suave de la variedad algebraica ${\mathcal M}$ es una variedad
simpl\'{e}ctica con forma simpl\'{e}ctica
\begin{eqnarray*}
\omega(A,B)=-\int_Str\, (A\wedge B).
\end{eqnarray*}
En esta f\'{o}rmula los vectores tangentes a ${\mathcal M}$, $A$ y $B$ admiten 
 dos conexiones de $su(2)$ tangentes al espacio af\'{\i}n
de todas la conexiones de $su(2)$ en la superficie, y despues se
multiplican y el producto se integra sobre la superficie.

\begin{figure}[htbp]
\centering
\scalebox{.30}{\includegraphics{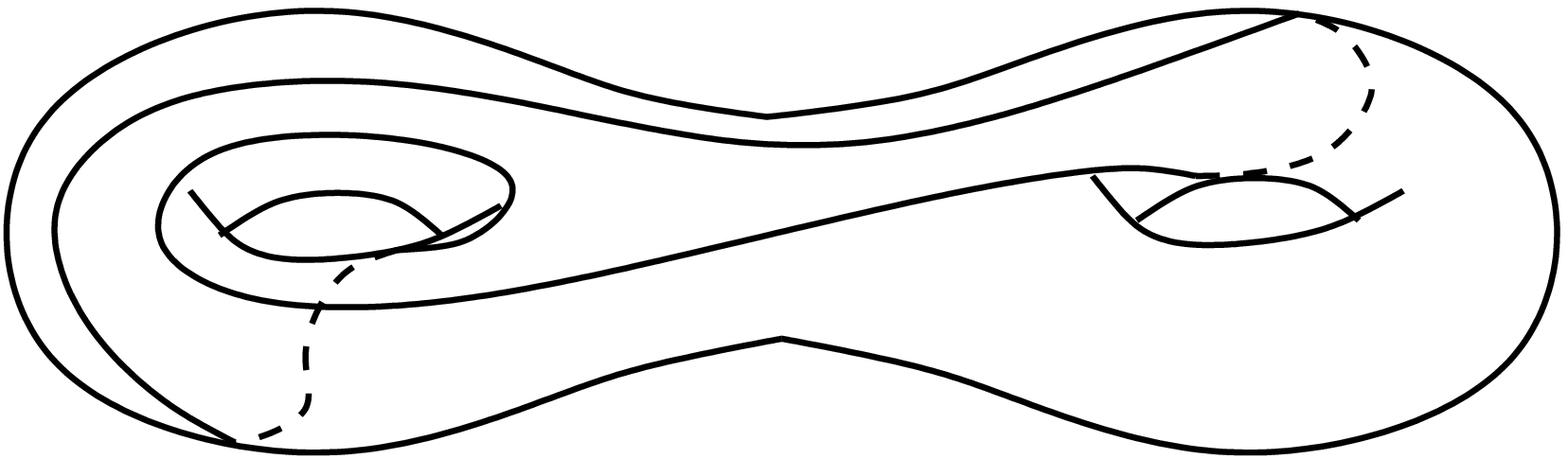}}

Figura 5.
\end{figure}

Los observables cl\'{a}sicos  son las trazas de holonom\'{\i}as 
sobre curvas cerradas (Figura 5) en la superficie:
\begin{eqnarray*}
 {I_{\gamma }:\, A\rightarrow tr\,hol_\gamma(A))}.
\end{eqnarray*}

W. Goldman ha mostrado \cite{goldman}
que el par\'{e}ntesis de Poisson 
(necesario para hacer mec\'{a}nica cl\'{a}sica en el espacio
${\mathcal M}$) se calcula con
la f\'{o}rmula 
\begin{eqnarray*}
{\{I_{\alpha},I_{\beta}\}= \frac{1}{2}
\sum_{x\in \alpha\cap \beta}
sgn(x)(I_{\alpha\beta^{-1}_x}-I_{\alpha\beta_x})}
\end{eqnarray*}
donde el signo del cruzamiento $x$ es ``$-$'' si $\alpha $ y $\beta$ se cruzan
como en la Figura 6, y ``$+$'' si $\beta $ tiene la direcci\'{o}n opuesta.
Aqui pensamos de las curvas $\alpha$, $\beta$, $\alpha\beta^{-1}_x$,
y $\alpha\beta_x$ como elementos del grupo fundamental de la superficie.

\begin{figure}[htbp]
\centering
\scalebox{.27}{\includegraphics{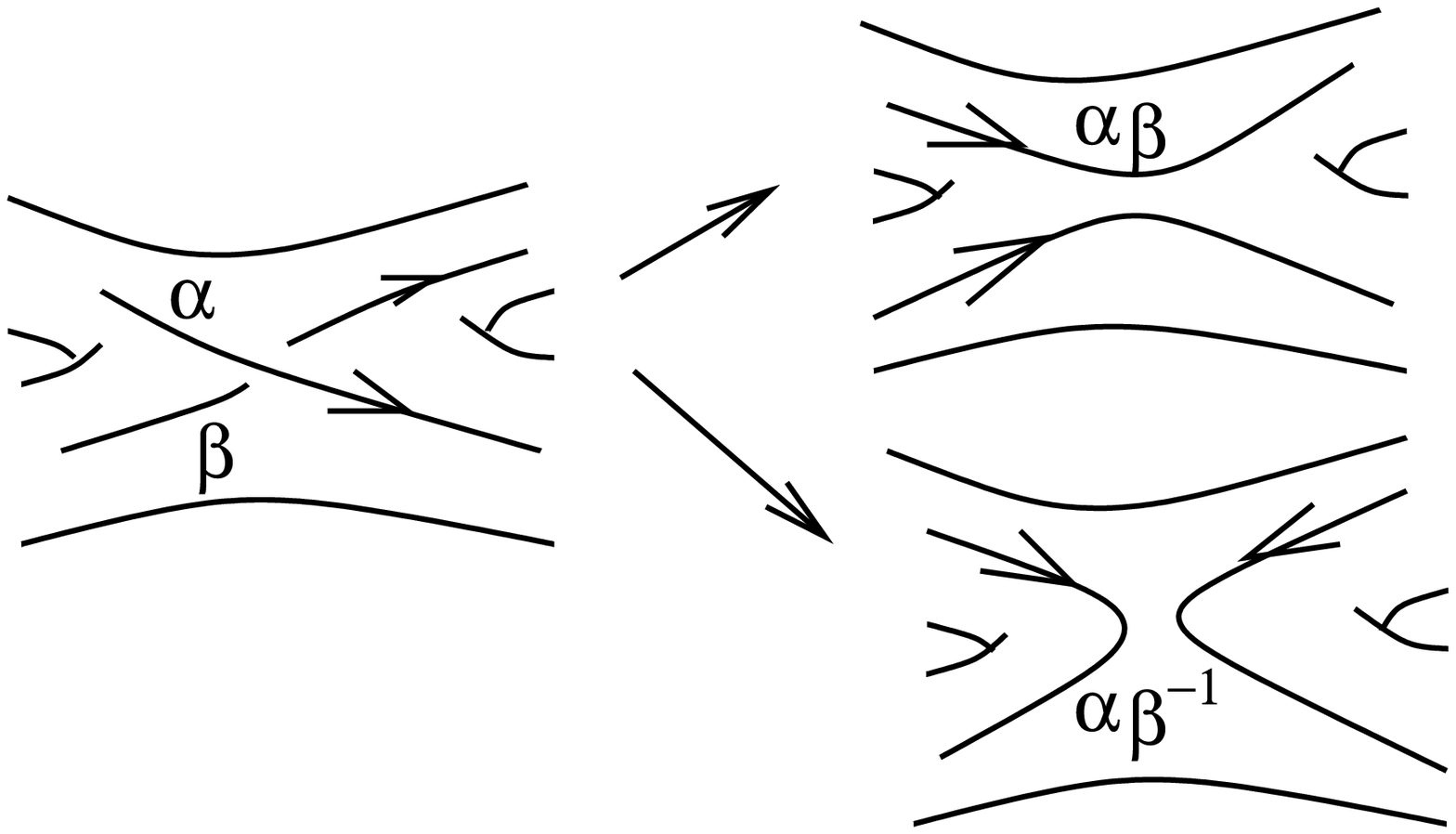}}

Figura 6.
\end{figure}

El espacio de Hilbert   ${\mathcal H}_r$ de la 
cuantizaci\'{o}n de ${\mathcal M}$  
consiste de las funciones theta generalizadas. 
\'{E}stas son las secciones
holom\'{o}rfas del haz en l\'{\i}neas de  Chern-Simons 
\begin{eqnarray*}
{\mathcal L}={\mathcal A}_F^s\times_{\Theta}{\mathbb C}.
\end{eqnarray*}
Aqui ${\mathcal A}_F^s$ representa el espacio de las
 conexiones planas irreducibles,
y el cociclo de Chern-Simons es 
\begin{eqnarray*}
(A,z)\sim (A^g,\Theta(A, g)z)
\end{eqnarray*} 
donde  $A^g=g^{-1}Ag+g^{-1}dg$ y 
$\Theta(A,g)=\exp(iL(\widetilde{A^g})-L(\widetilde{A}))$,
con $\widetilde A$  una extensi\'{o}n de la connexi\'{o}n $A$ en una
variedad de dimension tres limitada por la superficie y 
$\widetilde{A^g}$  una extensi\'{o}n de $A^g$ en la misma variedad
que es ``gauge'' equivalente a $\widetilde{A}$. 


El espacio de fase ${\mathcal M}$ es de volumen finito. Por este motivo,
el espacio de Hilbert asociado tiene dimensi\'{o}n finita. 
\'{E}sta es una consecuencia del principio de incertidumbre de Heisenberg,
porque cada part\'{\i}cula ocupa un volumen positivo en el espacio de fase,
y el numero maximo de particulas es igual a la dimensi\'{o}n del espacio de
Hilbert.    

La  dimensi\'{o}n de ${\mathcal H}_r$ se calcula con la f\'{o}rmula de 
Verlinde  (probada para un grupo de Lie
general por Faltings \cite{faltings}): 
\begin{eqnarray*}
{dim \, {\mathcal H}_r=r^{g-1}\sum_{j=1}^{r-1}
\left(2\sin^2\frac{j\pi}{r}\right)^{-(g-1)},}
\end{eqnarray*}
donde  $g$ es el g\'{e}nero de la superficie.

Por otro lado, los observables cu\'{a}ntic{o}s 
son muy misteriosos. Witten dice \cite{witten} que
la cuantizaci\'{o}n de $tr \, hol_\gamma(A)$ (donde  $\gamma$ 
es una curva cerrada en la
superficie) es el  operador integral con  n\'{u}cleo
\begin{eqnarray*}
<A_1|tr \, hol_\gamma(A)|A_2>=\int_{\mathcal M_{A_1,A_2}}
e^{2irL(A)}tr\, hol_\gamma(A){\mathcal D}A.
\end{eqnarray*}

Esta definici\'{o}n usa la integral de Feynman, y por tanto
no tiene un fundamento matem\'{a}tico.
Ya hay un \'{u}nico modelo de cuantizaci\'{o}n de los operadores,
que es parte de la teor\'{\i}a topol\'{o}gica 
de campos cu\'{a}nticos de Reshetikhin y Turaev \cite{RT} 
 basada en los grupos cu\'{a}nticos de Drinfeld \cite{drinfeld}. Esta
teor\'{\i}a nos da una
construcci\'{o}n combinatoria de los operadores. Presentamos
en lo qu\'{e} sigue las   id\'{e}as b\'{a}sicas de esta 
cuantizaci\'{o}n.
   
El grupo cu\'{a}ntico de $SU(2)$ tiene un numero finito de 
 representaciones  irreducibles  $V^1,V^2,\ldots,
V^{r-1}$, donde el \'{\i}ndice representa la dimensi\'{o}n. 
Usando la teoria de campos cu\'{a}nticos 
conformes \cite{segal} 
se pueden identificar las funciones theta generalizadas en la
superficie  con  coloraciones de las orillas del grafo que es el
 centro de la superficie con representaciones
irreducibles. 
La Figura 7 muestra el caso de la superficie de g\'{e}nero $2$.
Las dimensiones $m$, $n$, $p$ de las tres representaciones que
se encuentran a cada v\'{e}rtice satisfacen las condiciones 
cuantizadas de Clebsch-Gordon $${|m-n|+1\leq p\leq \min(m+n-1,2r-2-m-n)}.$$

\begin{figure}[htbp]
\centering
\scalebox{.35}{\includegraphics{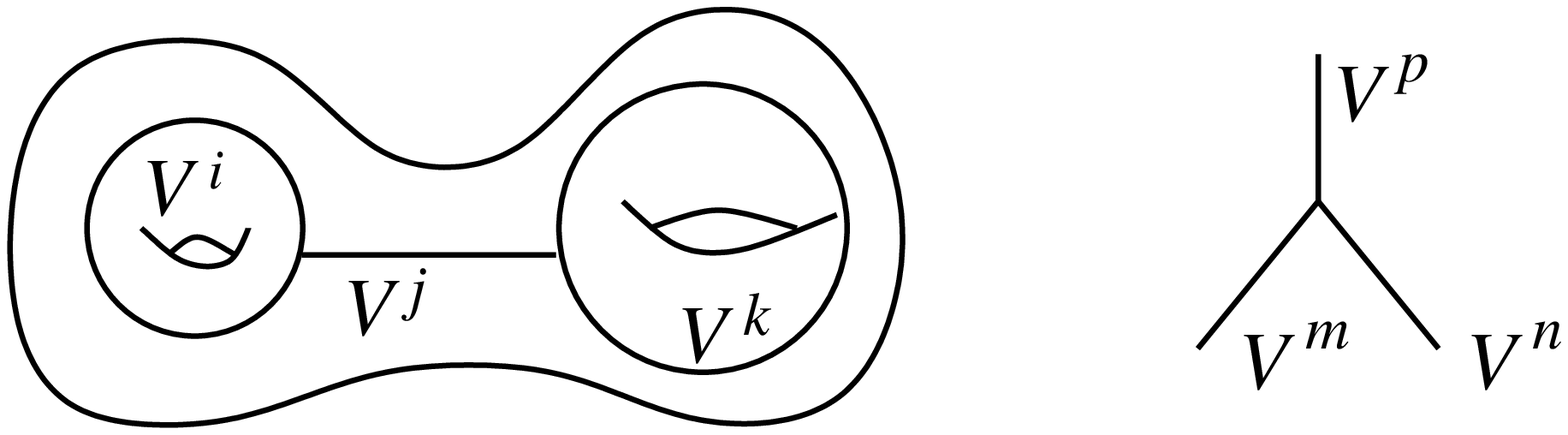}}

Figura 7.
\end{figure}

Por otro lado, las matrices de los operadores se calculan
 con diagramas de nudos y
de grafos con relaciones como la del polinomio de Jones \cite{RT}. 
Los calculos son muy complicados, y son basados en 
funciones hipergeom\'{e}tricas.

\begin{figure}[htbp]
\centering
\scalebox{.50}{\includegraphics{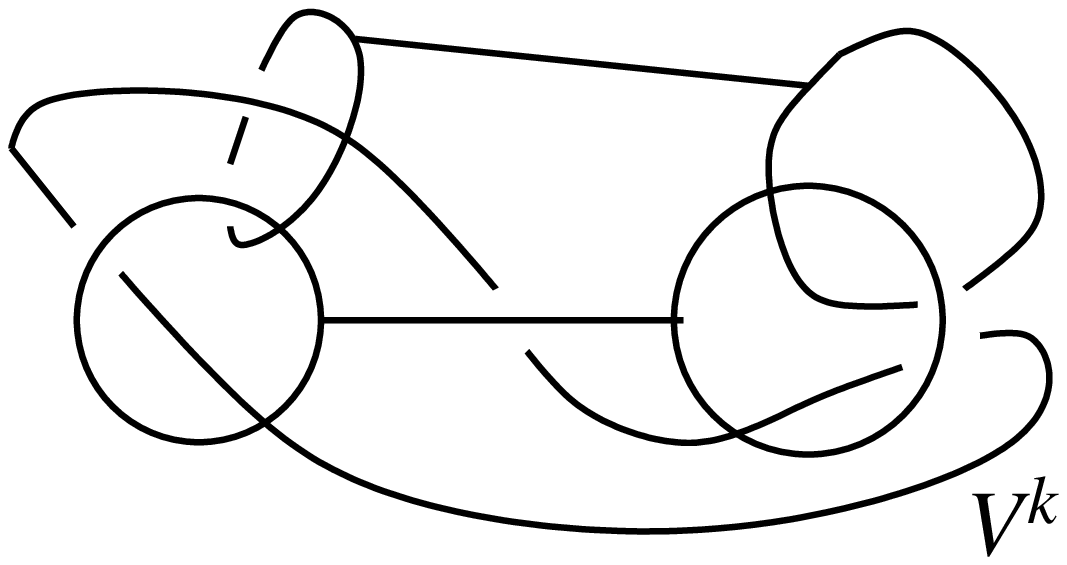}}

Figura 8.
\end{figure}

Por ejemplo el grafo con curva cerada colorada en la Figura 8 describe 
un operador para el caso de la superficie
de g\'{e}nero 2, y cada coloraci\'{o}n de las orillas 
de este grafo con representaciones
irreducibles de  $SU(2)$ cu\'{a}ntico produce una entrada de la matriz del
operador. 

Alexeev y Schomerus \cite{alexeevschomerus} han demonstrado que esta 
cuantizaci\'{o}n satisface la condici\'{o}n de Dirac 
\begin{eqnarray*}
op(\{I_{\alpha},I_{\beta}\})= \frac{1}{i\hbar}
[op (I_{\alpha}), op(I_{\beta})]+O(\hbar)
\end{eqnarray*}
para el par\'{e}ntesis de Poisson de Goldman. 
Pero nos interesa entender esta cuantizaci\'{o}n desde un punto
de vista geom\'{e}trico, con mas intuici\'{o}n f\'{\i}sica.
El \'{u}nico caso bien entendido es el del toro (por que el caso de la
esfera es trivial).

\section{El caso particular del toro}

En el caso del toro, el espacio de fase ${\mathcal M}$ se puede
calcular directamente y es muy sencillo.  Este espacio se
llama la {``funda de almohada''} y se obtiene factorizando
el plano complejo por las transformaciones $z\rightarrow z+m+in$,
$m,n\in {\mathbb Z}$ y $z\rightarrow -z$ (Figura 9).

\begin{figure}[htbp]
\centering
\scalebox{.40}{\includegraphics{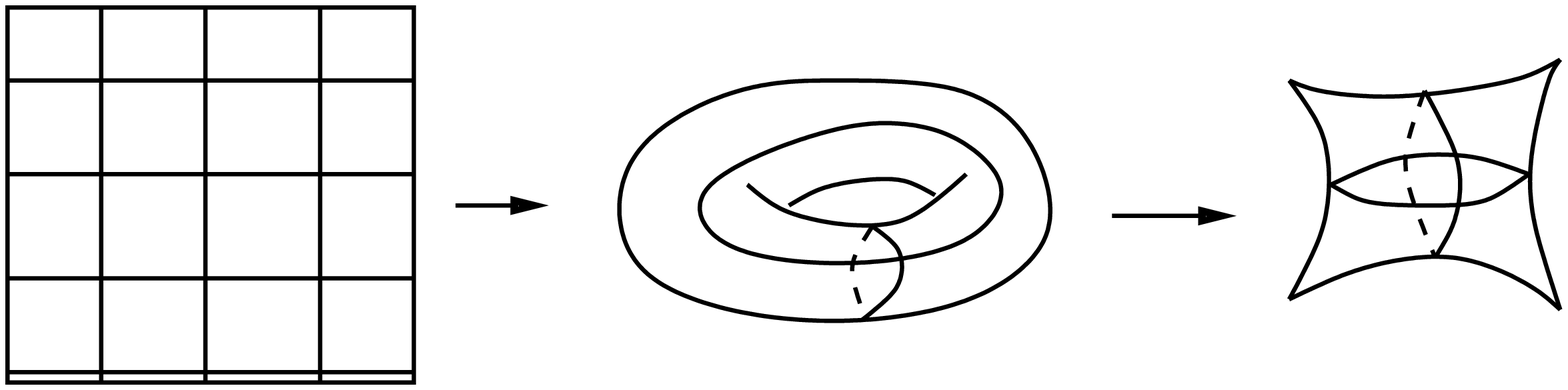}}

Figura 9.
\end{figure}

El autor de este art\'{\i}culo en su trabajo con Frohman \cite{FG} ha 
descubierto la f\'{o}rmula siguiente para  multiplicar los 
operadores $C(p,q)$ asociados a los observables clasicos $2\cos 2\pi(px+qy)$
en la funda de almohada:
\begin{eqnarray*}
C(m,n)C(p,q)=e^{i\pi(mq-np)/r}C(m+p,n+q)+e^{-i\pi(mq-np)/r}C(m-p,n-q).
\end{eqnarray*}
Desde el punto de vista de la teor\'{\i}a de Reshetikhin y Turaev,
el operador $C(p,q)$ es la diferencia de los operadores definidos por
 la curva de 
inclinaci\'{o}n $p/q$ en el toro colorada por $V^{d+1}$ y $V^{d-1}$, donde
$d$ es el m\'{a}ximo com\'{u}n divisor de $p$ y $q$.  

Esta f\'{o}rmula  
esta relacionada con el principio de incertidumbre de  Heisenberg. La
forma mas reconocida del principio de incertidumbre  es 
$\Delta p\Delta q\geq \hbar$, pero
su  forma b\'{a}sica es $op(p)op(q)-op(q)op(p)=\hbar\, op(1)$. Si denotamos
${U}$ y ${V}$ los 
operadores de evoluci\'{o}n de la posici\'{o}n
y del momento
\begin{eqnarray*}
U=e^{i\,op(q)}, V=e^{i\,op(p)},
\end{eqnarray*}
entonces el principio de incertidumbre se convierte en 
\begin{eqnarray*}
VU=e^{2\pi i\hbar}UV.
\end{eqnarray*}
El \'{a}lgebra von Neumann generada por $U, V, U^{-1}$ y $V^{-1}$ se
llama el {toro no conmutativo} y es una deformaci\'{o}n del \'{a}lgebra
$C^{\infty}(S^1\times S^1)$ de las funciones en el toro. El toro no conmutativo
 ocurri\'{o} por primera vez
en el trabajo de Rieffel \cite{Ri} y  juega un rol central
en la geometr\'{\i}a  no conmutativa de Connes \cite{connes}.

El resultado con Frohman demuestra que hay un morfismo del \'{a}lgebra 
de los operadores de la cuantizaci\'{o}n de ${\mathcal M}$ en el 
toro no conmutativo 
\begin{eqnarray*}
op(2\cos 2\pi (px+qy))\longrightarrow
e^{-\pi i\hbar pq}(U^pV^q+U^{-p}V^{-q}).
\end{eqnarray*}
Recuerde que $\hbar=\frac{1}{2r}$. 

Este resultado tiene una explicaci\'{o}n mas profunda. En su trabajo con
Uribe \cite{GU} el autor de este art\'{\i}culo ha demonstrado que 
la cuantizaci\'{o}n de ${\mathcal M}$ con grupos cuanticos  es 
isomorfa a la cuantizaci\'{o}n de este espacio con el procedimiento de
H. Weyl (la que usa la trasformada de Fourier).

Para describir este isomorphismo utilizamos la forma compleja de la 
cuantizaci\'{o}n de Weyl.
Las funciones de onda que son las coloraciones del
centro del toro con $V^j$, $j=1,2,\cdots, r-1$, como mostrado en
la Figura 10, 
coresponden a las funciones theta 
\begin{eqnarray*}
{\zeta_j(z)}=
\sqrt[4]{r}e^{-\frac{\pi j^2}{2r}}({\theta_j(z)}-
{\theta_{-j}(z)}), \quad j=1,2,\ldots, r-1,
\end{eqnarray*}
donde 
\begin{eqnarray*}
\theta_j(z)=\sum_{n=-\infty}^{\infty}e^{-\pi(2rn^2+2jn)+2\pi iz(j+2rn)}.
\end{eqnarray*}

\begin{figure}[htbp]
\centering
\scalebox{.30}{\includegraphics{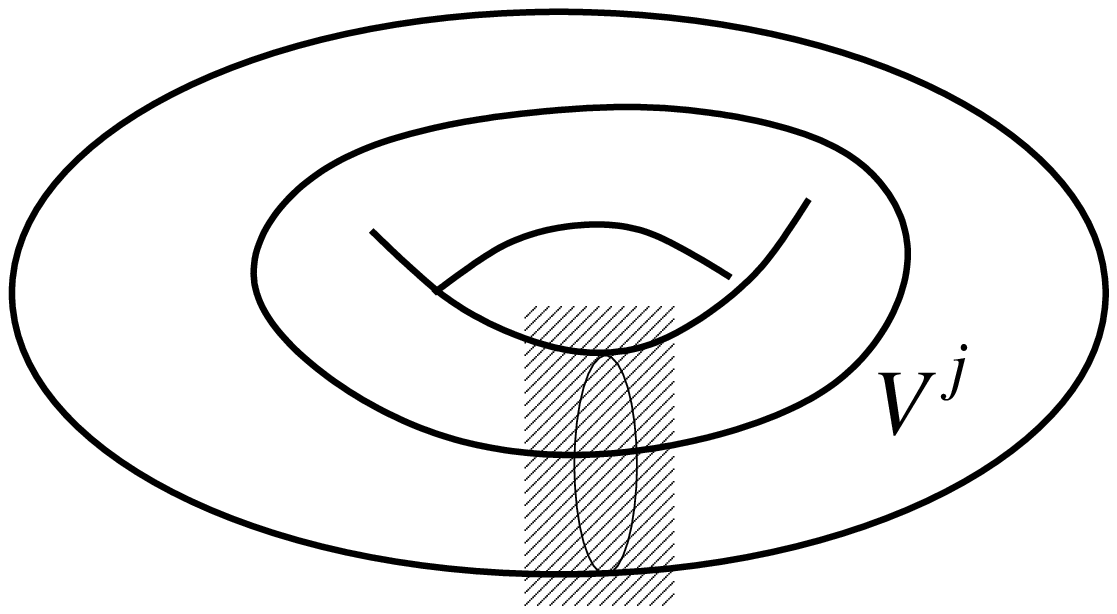}}

Figura 10.
\end{figure}

En la cuantizaci\'{o}n con grupos cuanticos, el operador 
asociado a la traza en la representaci\'{o}n en dimension $2$
de la holonom\'{\i}a
por una  curva con inclinaci\'{o}n 
 $p/q$ en el toro es dado por la f\'{o}rmula
\begin{eqnarray*}
\zeta_j\rightarrow e^{-\frac{i\pi}{2r}pq}(e^{\frac{i\pi}{r}jq}
\zeta_{j-p}+e^{-\frac{i\pi}{r}jq}\zeta_{j+p}).
\end{eqnarray*} 
%
En la  cuantizaci\'{o}n de Weyl
 el mismo  operador es el operador de Toeplitz con s\'{\i}mbolo
{$e^{-\Delta \hbar/4}\cos (px+qy)$}:
\begin{eqnarray*}
f\rightarrow e^{-\Delta \hbar/4}\cos 2\pi(px+qy)f\rightarrow 
\Pi_{{\mathcal H}_r} e^{-\Delta \hbar/4}\cos 2\pi(px+qy)f,
\end{eqnarray*}
donde $\Delta = \frac{1}{2\pi}(\partial^2/\partial x^2+\partial^2/\partial y^2)$ y $\Pi_{{\mathcal H}_r}$ es la proyecci\'{o}n en el espacio de
Hilbert ${\mathcal H}_r$.

\section{Una posible aplicaci\'{o}n}

Para motivar m\'{a}s el estudio de la teor\'{\i}a f\'{\i}sica del
polinomio de Jones, concluimos este art\'{\i}culo con una
aplicaci\'{o}n posible de esta teor\'{\i}a a la investigaci\'{o}n
del efecto Hall  cu\'{a}ntico fraccionario.

El efecto Hall cl\'{a}sico fue observado en 1879 de E. Hall. Este
efecto ocurre cuando una corriente el\'{e}ctrica pasa por un material
 en presencia de un campo
magn\'{e}tico (ver Figura 11). En este caso los electones tienden a irse hacia
un lado. La teor\'{\i}a cl\'{a}sica pronostica que la intensidad
de este efecto es proporcional a la intensidad del campo
magn\'{e}tico.  Pero en 1979, von Kitzling
descubri\'{o} que, para un campo magn\'{e}tico muy fuerte (10 T) 
y temperatura baja (50 mK) el efecto Hall est\'{a} cuantizado. Precisamente,
si $R_H=V_H/I$ denota la resistencia Hall, entonces 
$1/R_H$ esta cuantizado en unidades de
$e^2/h$. En los 80's, Tsui y Stormer descubrieron adem\'{a}s que
exist\'{\i}a una resistencia de Hall para fracciones de $e^2/h$, lo que
se llama el efecto Hall cu\'{a}ntico fraccionario. Este efecto
se observa cuando el flujo de electrones est\'{a} ligado a una
superficie, como por ejemplo en la juntura entre un cristal de
galio-ars\'{e}nico y un cristal de aluminio-galio-ars\'{e}nico.

 \begin{figure}[htbp]
\centering
\scalebox{.40}{\includegraphics{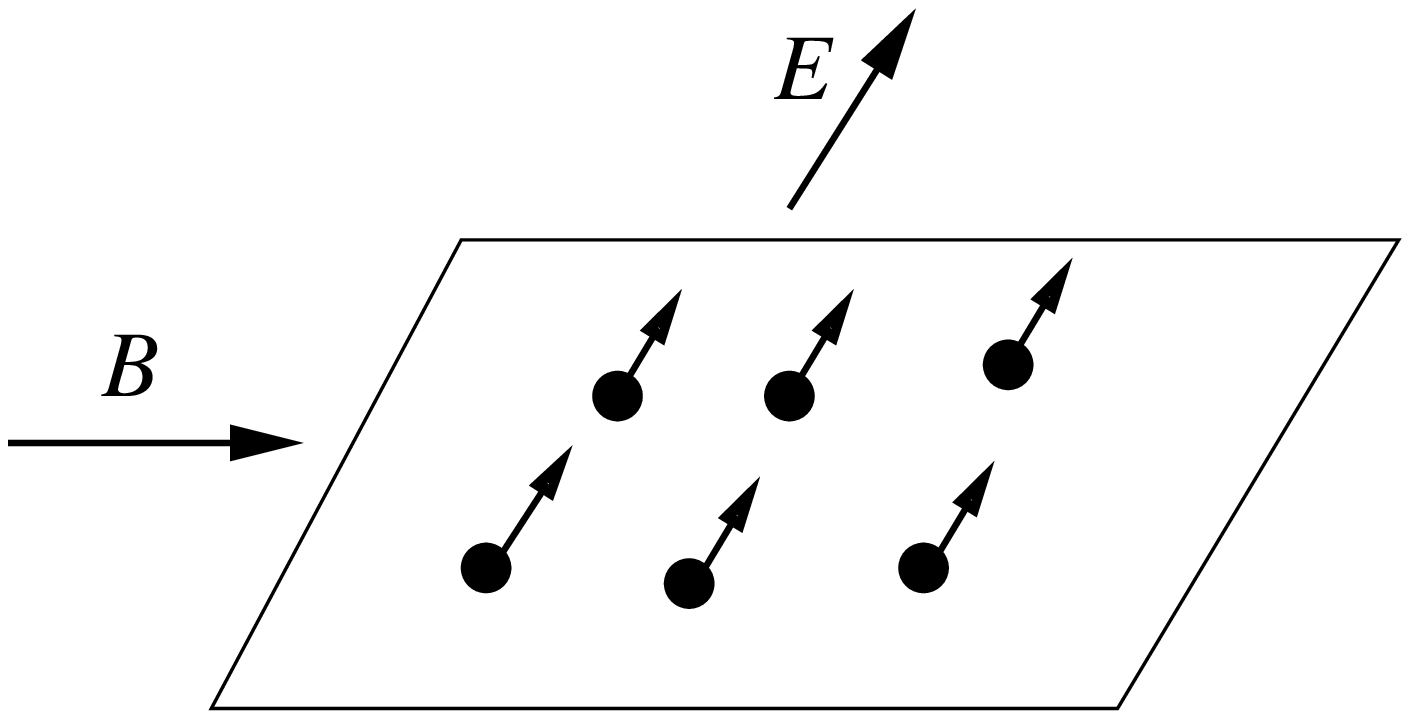}}

Figura 11.
\end{figure}





R. Laughlin  gan\'{o} el Premio Nobel para  explicar las fracciones,
llamadas niveles, ${\frac{1}{2m+1}}$, $m=1,2,3,\ldots$.
Su teor\'{\i}a pronostica la existencia de cargas 
el\'{e}ctricas fraccionales, con consecuencia la 
 existencia de la estad\'{\i}stica fraccionaria de part\'{\i}culas.
Pues adem\'{a}s de los bosones, cuyas funciones de onda
quedan invariantes por la permutaci\'{o}n de part\'{\i}culas, y de los
fermiones, cuyas funciones de onda cambian de signo por la
permutaci\'{o}n de part\'{\i}culas, hay otras que se
llaman cuasi-part\'{\i}culas, cuyas funciones de onda se multiplican
con un n\'{u}mero complejo cuando las cuasi-part\'{\i}culas se
permutan (ver Figura 12).    

\begin{figure}[htbp]
\centering
\scalebox{.40}{\includegraphics{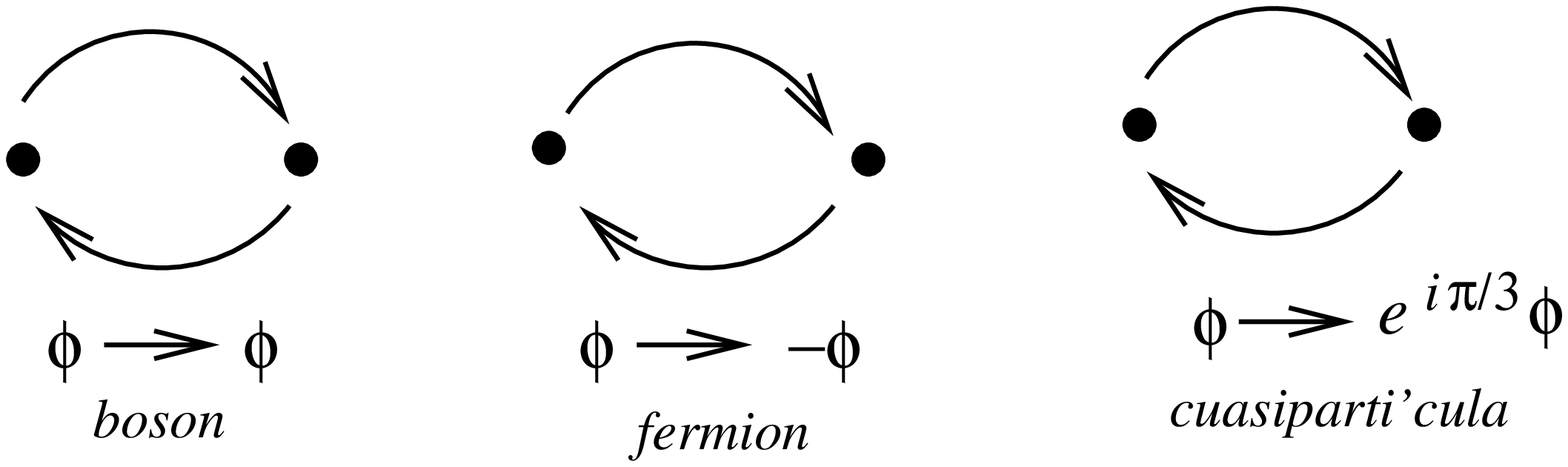}}

Figura 12.
\end{figure}

Moore y Read \cite{mooreread} conjecturaron que 
algunos niveles (e.g.~$\frac{5}{2}$) se pueden modelar
con la teor\'{\i}a de campo cu\'{a}ntico con lagrangiano de Chern-Simons,
la teor\'{\i}a f\'{\i}sica del polinomio de Jones.
El modelo  de Moore y Read pronostica
la existencia de  estad\'{\i}stica no Abeliana.
En este caso un sistema de part\'{\i}culas cuyas trayectorias forman
una trenza (Figura 13) produce una representaci\'{o}n del grupo de trenzas en
un grupo de matrices, y esta  es la representaci\'{o}n de Jones \cite{jones2}
relacionada a su polinomio. 
Suponiendo la conjectura  verdadera, Freedman, Nayak, Shtengel, Walker, Wang 
dedujeron \cite{fnsww} que se puede
construir un  computador cu\'{a}ntico universal que utiliza
el efecto Hall cu\'{a}ntico fraccionario.

\begin{figure}[htbp]
\centering
\scalebox{.40}{\includegraphics{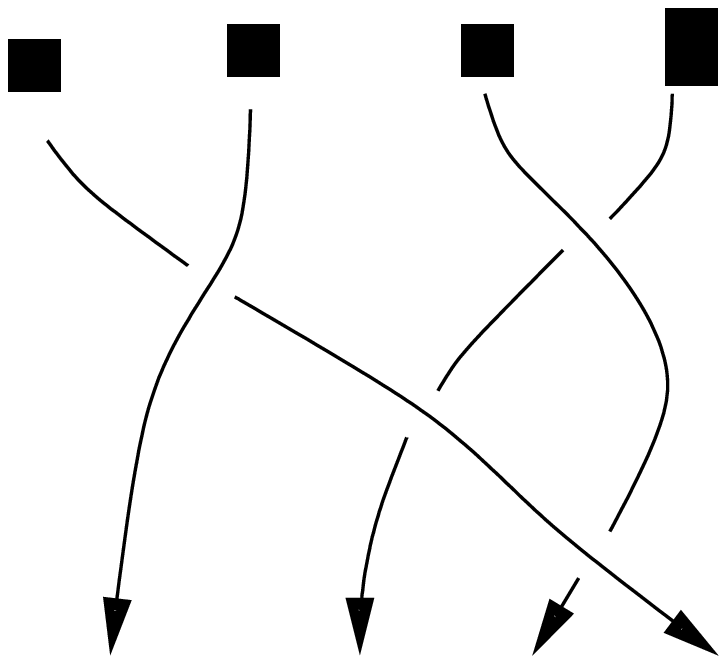}}

Figura 13.
\end{figure}




\end{document}